\documentclass[11pt]{scrartcl}
\usepackage{amsmath}
\usepackage{amssymb}
\usepackage{latexsym}
\usepackage{array}
\usepackage{yfonts}
\usepackage{fancyhdr}
\usepackage{epic, eepic}
\usepackage{graphicx,graphics}
\usepackage{color}
\usepackage{enumerate}
\DeclareMathAlphabet{\mathpzc}{OT1}{pzc}{m}{it}

\DeclareMathOperator{\Id}{Id}

\newtheorem{theorem}{Theorem}[section]
\newtheorem{defi}[theorem]{Definition}
\newtheorem{lemma}[theorem]{Lemma}

\newcommand{\R}{\mathbb{R}}

\newcommand{\N}{\mathbb{N}}

\newcommand{\beg}{\begin{equation}}
\newcommand{\en}{\end{equation}}
\newcommand{\beqn}{\begin{eqnarray}}
\newcommand{\eeqn}{\end{eqnarray}}
\newcommand{\beqnn}{\begin{eqnarray*}}
\newcommand{\eeqnn}{\end{eqnarray*}}

\numberwithin{equation}{section}

\newcommand{\tu}{\tilde{u}}

\author{}
\title{}
\date{}
\begin{document}
\maketitle
\vspace{-3.1cm}
\begin{center}{\LARGE $C^1$-regularity for local graph representations of immersions}\end{center}
\vspace{7mm}
\begin{center} Patrick Breuning \footnote{P.\ Breuning was supported by the DFG-Forschergruppe
\emph{Nonlinear Partial Differential Equations: Theoretical and Numerical Analysis}.
The contents of this paper were part of the author's dissertation, which was written at
Universit\"{a}t Freiburg, Germany.}
\\ Institut f\"{u}r Mathematik der Goethe Universit\"{a}t Frankfurt am Main \\ Robert-Mayer-Stra{\ss}e 10,
D-60325 Frankfurt am Main, Germany \\email: breuning@math.uni-frankfurt.de \end{center} \vspace{1mm}
\begin{abstract}
\begin{center} {\textbf{Abstract}} \end{center}
\noindent We consider immersions admitting uniform graph representations over the affine tangent space
over a ball of fixed radius $r>0$. We show that for sufficiently small $C^0$-norm of the graph functions,
each graph function is smooth with small $C^1$-norm.
\end{abstract}
\begin{section}{Introduction}
An immersion into $\R^n$ is a differentiable function $f:M\rightarrow \R^n$ defined on a differentiable
manifold $M^m$, such that for each
$q\in M$ the mapping $f_\ast|T_qM$ is injective.
A simple consequence of the implicit function theorem says that any immersion
can locally be written as the graph of a function $u:B_r\rightarrow\R^k$ over the affine tangent space. Moreover,
for a given $\lambda>0$ we can choose $r>0$ small enough such that $\|Du\|_{C^0(B_r)}\leq \lambda$. If
this is possible at any point of the immersion with the same radius $r$, we call $f$ an $(r,\lambda)$-immersion. \\ \\
This concept is used in various geometric contexts; as an example and as motivation we consider the following
compactness theorem proved by J.\ Langer \cite{langer}: Let $f^i:\Sigma^i\rightarrow \R^3$ be a sequence
of immersed surfaces with uniformly $L^p$-bounded second fundamental form, $p>2$, and uniformly bounded area.
Then, after passing to a subsequence, there are a limit immersion $f:\Sigma\rightarrow \R^3$ and
diffeomorphisms $\phi^i:\Sigma\rightarrow \Sigma^i$, such that $f^i\circ \phi^i$ converges in the $C^1$-topology to $f$.
The result can be generalized to higher dimensions and codimensions; see
\cite{breuning3}, \cite{corlette}, \cite{delladio}.
For proving the statement, one uses the Sobolev embedding and shows that a uniform $L^p$-bound
for the second fundamental form with $p$ greater
than the dimension implies that for any $\lambda>0$ there is an $r>0$ such that every immersion is an
$(r,\lambda)$-immersion. \\ \\
This conclusion plays an important role in the proof of the compactness theorem and is just one
example of a fundamental principle frequently used in geometric analysis and related fields: For a given
global object, that is a manifold embedded or immersed in $\R^n$ --- usually of some specific geometric type,
for example a minimal surface --- one investigates the local graph representations in order to derive
further characteristics of the given object.
For that one uses the global geometric information and derives specific properties satisfied by
each of the graph functions, for example bounds for specific norms, or particular partial
differential equations to be satisfied. For each of the graph functions, it is then possible to apply all the
well-known results from real analysis like embedding theorems or regularity theory. \\ \\
In this paper, we like to take a slightly different point of view. Instead of deriving
special kinds of graph representations from specific geometrical settings, we shall take immersions
with specific graph representations as our starting point. More precisely, our concept is the
following: We consider an immersion
and assume that it can be represented at any point over a ball of
fixed radius $r>0$ as the graph of
a function $u$ satisfying some specific properties; now, loosely speaking,
we claim that each of the graph functions satisfies much better properties than one would anticipate
from the ordinary rules of analysis. \\ \\
In fact, there is a huge difference between a single graph and a graph coming from an \nopagebreak
immersion in the way described above. In the latter case, we know as an additional information
that such a graph representation is possible at \emph{any} point of the immersion.
In particular, two graphs that are close to each other have overlapping parts and each of the graphs
satisfies specific properties, such as a bounded norm. Hence
all graphs having one point in common depend on each other. This can be seen as a combinatorial
restriction and allows much stronger results than one would expect using only the
given properties of each single graph. \\ \\
Let us first generalize the concept of immersions with bounded norm
$\|Du\|_{C^0(B_r)}\leq \lambda$ for the graph functions $u$
to immersions satisfying only a weaker bound. Again we
consider $C^1$-immersions with graph representations $u:B_r\rightarrow \R^k$ over
the affine tangent space, but this time we only assume that
$\|u\|_{C^0(B_r)}\leq r\lambda$. If such a representation is possible at every point for fixed
$\lambda$ and $r$, we say that $f$ is a $C^0$-$(r,\lambda)$-immersion. The factor $r$ on the right hand
side is necessary for scale-invariance. A graph function of a $C^0$-$(r,\lambda)$-immersion
does not need to be differentiable; this explains the notation that we use for this kind of immersion.
For the precise definitions and further details the reader is referred to Section 2.\\ \\
Of course it is completely impossible to derive Lipschitz estimates for a single
function satisfying only a $C^0$-bound,
even if the function is known to be smooth or if the $C^0$-norm
is particularly small. However, as we have claimed above, graph functions coming from
immersions in the described way have much better properties than a single function.
Denoting by $m$ the dimension of the manifold on which the
immersion is defined, we obtain the following theorem: \\
\begin{theorem} \label{embedding} (Embedding theorem for
$C^{0}$-$(r,\lambda)$-immersions) \\ \\ For every $m \in \N$
there is a $\Lambda=\Lambda(m)>0$, such that every
$C^{0}$-$(r,\lambda)$-immersion with $\lambda \leq \Lambda$ is also an
$(r,\frac{\lambda}{\Lambda})$-immersion. \\ \\
The constant $\Lambda$ can be given explicitly by\,
$\Lambda(m):=10^{-5}m^{-2}$.
\end{theorem} \vspace{6mm}
Hence a sufficiently small $C^0$-norm implies that each graph function is
smooth with small $C^0$-norm of $Du$, that is with small Lipschitz constant. Equivalently,
we can say that the space of $C^0$-$(r,\lambda)$-immersions embeds into the
space of $(r,\frac{\lambda}{\Lambda})$-immersions.
The statement
is true in arbitrary codimension and also for noncompact manifolds.\\ \\
As here we are assuming only a small $C^0$-norm, we obtain all at once
whole classes of new embedding theorems --- provided the functions
come from graph representations as described above.
For example one can think of the case of H\"{o}lder continuous
$C^{0,\alpha}$-graphs or the Sobolev border case of
$W^{2,m}$-graphs in dimension $m$. \\ \\
The question arises, whether the result will still be true, if we assume graph representations
not over the affine tangent space, but over other appropriately chosen $m$-spaces. In the appendix
we will show that this is not the case.
\\ \\ \\[-3mm]
\emph{Acknowledgement:} I would like to thank my advisor Ernst Kuwert for his support. Moreover
I would like to thank Manuel Breuning for proofreading
my dissertation \cite{breuning2}, where the result of this paper was established first.
\end{section}
\begin{section}{Notation and definitions}
We begin with some general notations: For $n=m+k$ let $G_{n,m}$ denote the Grassmannian
of (non-oriented) $m$-dimensional subspaces of $\R^n$. Unless stated otherwise let
$B_{\!\varrho}$ denote the open ball in $\R^m$
of radius $\varrho>0$ centered at the origin. \\ \\
Now let $M$ be an $m$-dimensional manifold without boundary and $f:M\rightarrow\R^n$
a $C^1$-immersion. Let $q\in M$ and let $T_qM$ be the tangent space at $q$. Identifying vectors
$X\in T_qM$ with $f_\ast X\in T_{f(q)}\R^n$, we may consider $T_qM$ as an $m$-dimensional
subspace of $\R^n$. In this manner
we define the tangent map
\\ \parbox{14.5cm}{\beqnn \tau_f:M&\rightarrow& G_{n,m}, \hspace{0cm}\\
q&\mapsto& T_qM.  \eeqnn}   \hfill  \parbox{8mm}{\beqn \label{notiontangentspace} \eeqn}
\\ \\
\textbf{The notion of an \boldmath$(r,\lambda)$\unboldmath-immersion:} \\ \\
We call a mapping $A:\R^n\rightarrow \R^n$ a \emph{Euclidean isometry}, if
there is a rotation $R\in \mathbb{SO}(n)$ and a translation $T\in\R^n$, such that
$A(x)=Rx+T$ for all $x\in \R^n$. \\ \\
\label{defeuclisom}For a given point $q \in M$ let $A_q: \mathbb{R}^n\rightarrow
\mathbb{R}^n$ be a Euclidean isometry,
which maps the origin to
$f(q)$, and the subspace $\mathbb{R}^m\times\{0\}\subset
\mathbb{R}^m \times \mathbb{R}^k$ onto $f(q)+\tau_f(q)$.
Let $\pi:\R^n\rightarrow \R^m$ be the standard projection onto the
first $m$ coordinates. \\ \\
Finally let $U_{r,q}\subset M$ be the $q$-component
of the set \label{qcomponentof} $(\pi \circ A_q^{-1} \circ f)^{-1}(B_r)$.
Although the isometry $A_q$ is not uniquely determined, the set
$U_{r,q}$ does not depend on the choice of $A_q$.\\ \\
We come to the central definition (as first defined in \cite{langer}):
\begin{defi} \label{definition1rlambda} An immersion $f$ is called an
$(r,\lambda)$-immersion, if for each point $q\in M$ the set
  $A_q^{-1}\circ f(U_{r,q})$ is the graph of a differentiable function
  $u:B_r\rightarrow \R^k$ with $\|Du\|_{C^{0}(B_{r})} \leq \lambda$. \end{defi}
\vspace{3mm}
\setlength{\unitlength}{0.7cm}
\hspace{0cm}
\begin{picture}(12,12) \put(4.7,2){
\qbezier(9,5)(7.5,8.5)(12,8.5) \qbezier(12,8.5)(16,8.5)(13.5,5)
\qbezier(13.5,5)(12.2,3.2)(9,5) \qbezier(9,5)(6.5,6.5)(7.5,3.5)
\qbezier(7.5,3.5)(8.5,1)(9.5,1.5) \qbezier(9.5,1.5)(10.4,1.9)(9,5)

\qbezier(9,-7)(6.015855,-4.634655)(9.09463,-2.969985)
\qbezier(9.09463,-2.969985)(13.49288,-0.591885)(12.958426,-4.85971)
\qbezier(12.958426,-4.85971)(12.670997,-7.061386)(9,-7)
\qbezier(9,-7)(6.087445,-6.869575)(8.393955,-9.0329)
\qbezier(8.393955,-9.0329)(10.462655,-10.75641)(11.104495,-9.84097)
\qbezier(11.104495,-9.84097)(11.705932,-9.06105)(9,-7)

\put(9,-7){\line(1,0){4.5}} \put(9,-7){\line(-1,0){2.5}}

\qbezier(-2.25,5)(-2.25,6.1391)(-1.4445,6.9445)
\qbezier(-1.4445,6.9445)(-0.6391,7.75)(0.5,7.75)
\qbezier(0.5,7.75)(1.6391,7.75)(2.4445,6.9445)
\qbezier(2.4445,6.9445)(3.25,6.1391)(3.25,5)
\qbezier(3.25,5)(3.25,3.8609)(2.4445,3.0555)
\qbezier(2.4445,3.0555)(1.6391,2.25)(0.5,2.25)
\qbezier(0.5,2.25)(-0.6391,2.25)(-1.4445,3.0555)
\qbezier(-1.4445,3.0555)(-2.25,3.8609)(-2.25,5)

\linethickness{0.7mm}
\qbezier(-2.25,5.3)(-2.2,6)(-1.75,6.6)
\qbezier(0.5,7.75)(1.55,7.75)(2.4,6.95)
\qbezier(3.23,5.2)(3.3,4)(2.65,3.3)
\qbezier(1,2.3)(-0.2,2.095)(-1,2.7)

\thinlines \put(-2.5,-7){\line(1,0){7}}

\put(9,-7.1){\tiny $|$} \put(1,-7.1){\tiny $|$} \put(0,-7.15){(}
\put(1.9,-7.15){)} \put(8.92,-7.6){0} \put(1.9,-7.7){$B_r$}
\normalsize

\qbezier(3,8.8)(6,10.8)(9,8.8) \qbezier(14,2.5)(16,-0.5)(14,-3.5)
\qbezier(9,-10)(6,-12)(3,-10)
\qbezier(-2.5,2.5)(-4.5,-0.5)(-2.5,-3.5)

\put(9,8.8){\vector(3,-2){0.1}} \put(14,-3.5){\vector(-2,-3){0.1}}
\put(3,-10){\vector(-3,2){0.1}} \put(-2.5,-3.5){\vector(2,-3){0.1}}

\put(6,8.86){$f$} \put(13.8,-0.5){$A_q^{-1}$} \put(5.8,-10.7){$\pi$}
\put(-3,-0.5){$\pi\circ A_q^{-1}\circ f$}

\put(-2,7.8){$M^m$} \put(13,9){$\mathbb{R}^n=\mathbb{R}^m\times
\mathbb{R}^k$} \put(13.7,-7.2){$\mathbb{R}^m$}

\put(3.095,4.5){$\bullet$} \put(2.75,4.5){$q$}
\put(3,3.1){$U_{r,q}$}

\put(8.93,4.92){\tiny $\bullet$} \put(9.9,5.9){\vector(-1,-1){0.75}}
\put(9.75,6.15){$f(q)$}}
\end{picture}
\vspace{6.7cm} \\ \label{figureimmersion}
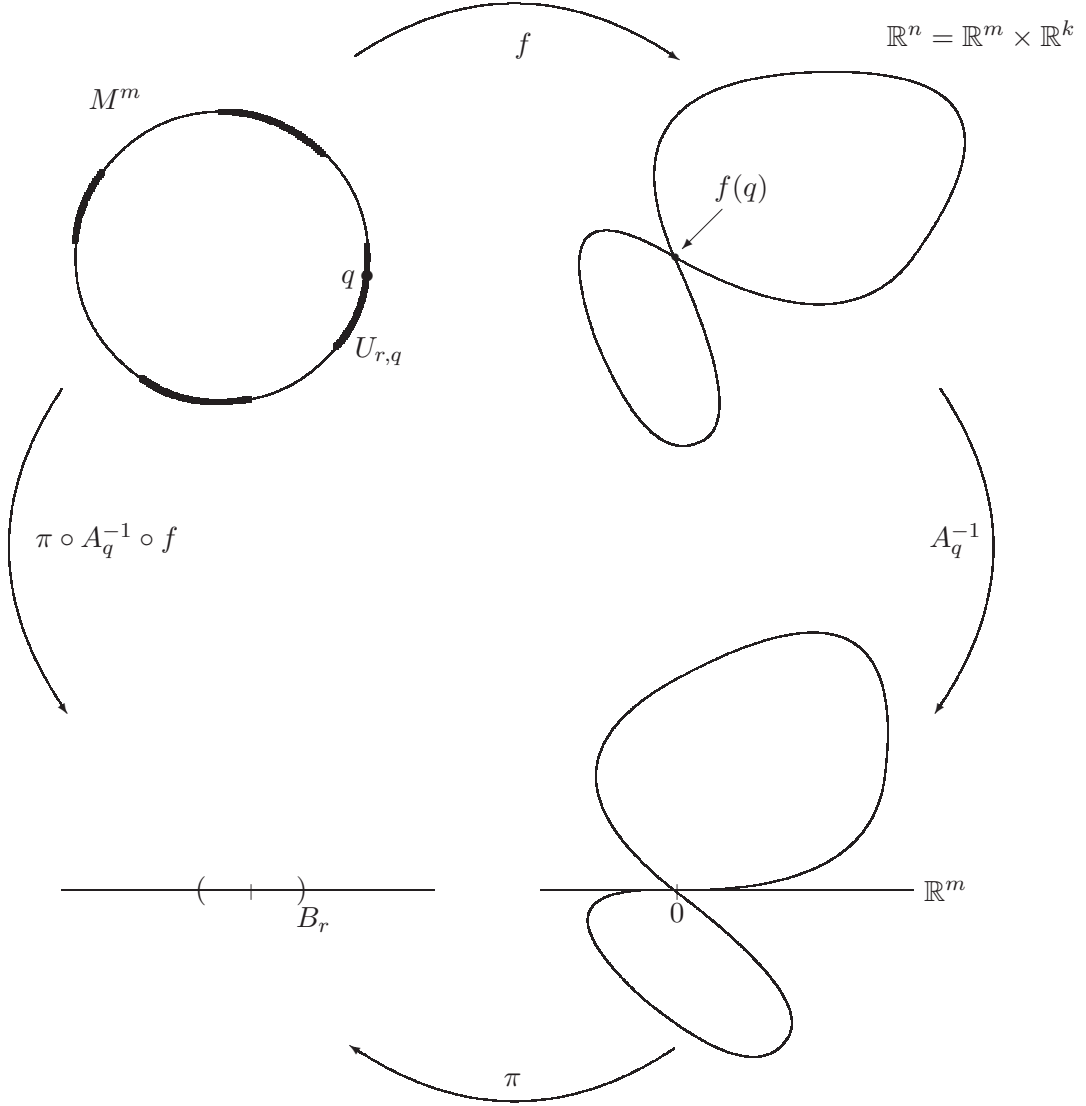
\noindent \textbf{Figure 2.1} \emph{Local representation as a graph.
The subset of $M$ drawn in bold lines
represents the pre-image $(\pi\circ A_q^{-1}\circ f)^{-1}(B_r)$.} \\ \\
Here, for any $x\in B_r$ we have $Du(x)\in \R^{k\times m}$. In order to define the
$C^0$-norm for $Du$, we have to fix a matrix norm for $Du(x)$. Of course all norms on
$\R^{k\times m}$ are equivalent, therefore our results are true for any norm
(possibly up to multiplication by some positive constant). Let us agree upon \vspace{-4mm}
\beqnn \|A\|\:=\:\Biggl(\hspace{0.5mm}\sum_{j=1}^m|a_j|^2\Biggr)^{\!\frac{1}{2}} \eeqnn
\vspace{-4mm} \\ for $A=(a_1,\ldots,a_m)\in \R^{k\times m}$. For this norm we have
$\|A\|_{\text{op}}\:\leq\: \|A\|$ for any $A\in \R^{k\times m}$
and the operator norm $\|\cdot\|_{\text{op}}$.
Hence the bound $\|Du\|_{C^0(B_r)}\leq \lambda$ directly implies that $u$
is $\lambda$-Lipschitz.
Moreover the norm $\|Du\|_{C^{0}(B_{r})}$ does not depend on the choice of the
isometry \nolinebreak$A_q$. \pagebreak \\
\textbf{The notion of a \boldmath$C^0$-$(r,\lambda)$\unboldmath-immersion:} \\ \\
Every $(r,\lambda)$-immersion admits a local representation as a graph of a differentiable function $u$
with $\|Du\|_{C^0(B_r)}\leq \lambda$.
This inequality corresponds to an estimate of the slope of the graph, i.e.\ to an estimate
of the Lipschitz constant of $u$. It is a natural generalization to consider immersions
with graph functions $u$, which satisfy only a bound for
some weaker norm. Any such definition should reasonably be scale-invariant
(i.e.\ if $f$ is an $(r,\lambda)$-immersion and $c>0$, then $cf$
is a $(cr,\lambda)$-immersion). \\ \\
Assuming only a bound for the $C^0$-norm yields the notion of a $C^0$-$(r,\lambda)$-immersion: \\[-3mm]
\begin{defi} \label{c0rlambda} An immersion $f$ is called a
$C^{0}$-$(r,\lambda)$-immersion, if for each point $q\in M$ the set
 $A_q^{-1}\circ f(U_{r,q})$ is the graph of a continuous function
 $u:B_r\rightarrow \R^k$ with $\|u\|_{C^{0}(B_{r})} \leq
r\lambda$.
\end{defi} \vspace{4mm}
It would not be sensible here to assume
$\|u\|_{C^{0}(B_{r})} \leq \lambda$, as the notion
of $C^0$-$(r,\lambda)$-immersions would not be scale-invariant then. For that reason we
require the bound $r\lambda$. \\ \\
Here we require $u$ only to be a continuous function. Note that the assumption on $f$ to be a smooth
immersion does \emph{not}\, imply that $u$ is differentiable. Surely the implicit function theorem
ensures a smooth graph representation over the tangent space. However this representation might
only be possible for radii less than $r$. Over the ball $B_r$ one might have a continuous
graph representation with a graph which gets vertical in a point. Hence smoothness of $f$
does not guarantee smoothness of $u$. \\
\vspace{0.3cm} \\ $\text{ }$
\hspace{0cm}
\begin{picture}(12,5.9)
\setlength{\unitlength}{0.45cm}
\qbezier(21.5,5)(21.5,8)(25.5,8)
\qbezier(21.5,5)(21.5,3)(23.5,3)
\qbezier(23.5,3)(24.5,3)(24.5,2)
\qbezier(24.5,2)(24.5,0)(26.5,0)
\qbezier(26.5,0)(29.5,0)(29.5,4)
\qbezier(29.5,4)(29.5,8)(25.5,8)

\qbezier(21.5,-5)(21.5,-2)(25.5,-2)
\qbezier(21.5,-5)(21.5,-7)(23.5,-7)
\qbezier(23.5,-7)(24.5,-7)(24.5,-8)
\qbezier(24.5,-8)(24.5,-10)(26.5,-10)
\qbezier(26.5,-10)(29.5,-10)(29.5,-6)
\qbezier(29.5,-6)(29.5,-2)(25.5,-2)

\put(23.4,-7.1){\tiny $\bullet$}
\put(23.5,-7){\line(1,0){1.8}}
\put(23.5,-7){\line(-1,0){1.8}}
\put(21.85,-7.1){$\underbrace{\text{ }}$}
\put(21.6,-7.1){\tiny $($}
\put(25.18,-7.1){\tiny $)$}
\put(22.5,-8.1){\small $r$}

\put(23.75,-7.75){\line(1,-4){0.82}}
\put(23.75,-7.75){\vector(-1,4){0.1}}
\put(24.6,-11.55){\small $f(q)$}

\put(23.5,-7){\circle{6}}
\put(23.5,-7){\circle{5.95}}
\put(23.5,-7){\circle{5.9}}
\put(23.5,-7){\circle{5.85}}
\put(23.5,-7){\circle{5.8}}
\put(23.5,-7){\circle{5.75}}
\put(23.5,-7){\circle{5.7}}
\put(23.5,-7){\circle{5.65}}
\put(23.5,-7){\circle{5.6}}
\put(23.5,-7){\circle{5.55}}

\qbezier(4.724,-5.362)(5.6522,-7)(8,-7)
\qbezier(8,-7)(9.82,-7)(9.82,-8.82)
\qbezier(9.82,-8.82)(9.82,-11.186)(11.458,-12.096)
\put(8,-7){\line(1,0){3.276}}
\put(8,-7){\line(-1,0){3.276}}

\put(5.9,-5.9){\line(3,4){0.9}}
\put(5.9,-5.9){\vector(-3,-4){0.1}}
\put(7.05,-4.6){\footnotesize $f(U_{r,q})$}

\put(9.52,-8.84){\line(-3,-1){1}}
\put(9.52,-8.84){\vector(3,1){0.1}}
\put(3.6,-9.35){\tiny graph gets vertical}
\put(4.1,-9.85){\tiny in one point}

\put(4.6,-7.15){\footnotesize $($}
\put(11.1,-7.15){\footnotesize $)$}

\qbezier(20.5,-9.5)(16.6,-12)(12.7,-9.5)
\put(12.7,-9.5){\vector(-3,2){0.1}}
\put(15.5,-10){zoom}

\put(8,4){\circle{8.5555}}
\put(3.5,7){$M$}

\qbezier(12.7,7)(16.6,9.5)(20.5,7)
\put(20.5,7){\vector(3,-2){0.1}}
\put(16.2,6.85){$f$}
\end{picture} \\
\vspace{5.31cm} \\
\label{graphnotsmooth} \noindent \textbf{Figure 2.2}
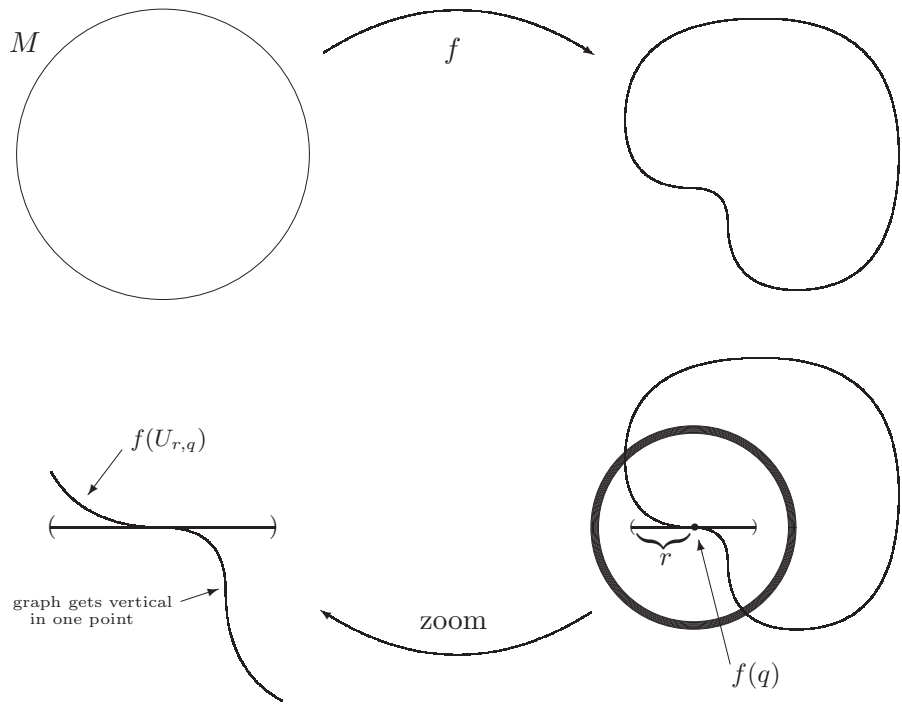
\emph{A simple example which shows how a graph function of a smooth
$C^0$-$(r,\lambda)$-immersion fails to be differentiable (here e.g. $\lambda=2$).} \pagebreak \\
\noindent Obviously every $(r,\lambda)$-immersion is also a $C^0$-$(r,\lambda)$-immersion.
Surprisingly, in some sense also the opposite is true: Every $C^0$-$(r,\lambda)$-immersion
is also an $(r,\frac{\lambda}{\Lambda})$-immersion if $\lambda\leq \Lambda=\Lambda(m)$. This
is precisely the statement of Theorem \ref{embedding}. Moreover, as we have seen above, a graph
function $u$ does not need to be smooth in the case of a $C^0$-$(r,\lambda)$-immersion; for that
reason we may interpret Theorem \ref{embedding} also as a higher regularity result.
\\ \\ \\
\textbf{Reformulation of Theorem \ref{embedding}:} \\ \\
Theorem \ref{embedding} is a statement for $C^0$-$(r,\lambda)$-immersions with fixed $r$ and $\lambda$.
We like to give an alternative formulation which holds for any immersion. \\ \\
For an immersion $f:M\rightarrow \R^n$
let $r_1(f,\lambda)\geq0$ be the maximal radius, such that for any $q\in M$ the set
$A_q^{-1}\circ f(U_{r,q})$ is the graph of a $C^1$-function $u:B_r\rightarrow \R^k$ with
$\|Du\|_{C^0(B_r)}\leq \lambda$. \\ \\
Similarly, let $r_0(f,\lambda)\geq 0$ be the maximal radius, such that for any $q\in M$ the set
$A_q^{-1}\circ f(U_{r,q})$ is the graph of a $C^0$-function with $\|u\|_{C^0(B_r)}\leq r\lambda$. \\ \\ \\
Obviously \beqnn
r_1(f,\lambda)\,\leq\: r_0(f,\lambda). \eeqnn \\[2mm]
With this notation Theorem \ref{embedding} reads as follows: \\[-1mm]
\begin{theorem}\label{reformembedding} (Reformulation of Theorem \ref{embedding}) \\ \\
For every $m \in \N$
there is a $\Lambda=\Lambda(m)>0$, such that for every
immersion $f:M^m\rightarrow \R^n$ and all $\lambda\leq\Lambda$ the inequality\,
$r_1(f,\lambda/\Lambda)\geq r_0(f,\lambda)$ holds. \\ \\
The constant $\Lambda$ can be given explicitly by\,
$\Lambda(m):=10^{-5}m^{-2}$. \end{theorem}
\end{section} \vspace{5mm}
\begin{section}{Preparations for the proof} \vspace{1mm}
The main step
of the proof is to compare the position of two tangent spaces at points on the surface that are
not too far from each other.
For that we have to find a sufficiently large set
$U\subset M$, such that $f(U)$ may be written over \emph{both} spaces
as graph with small $C^0$-norm respectively; this will be done in
Lemma \ref{inclusionc0}.
To compare the spaces with each other, we shall use
a finite number of comparison points on each space, constructed by means
of the immersion piece $f(U)$. A concrete estimate (in a slightly more general formulation)
is deduced in Lemma \nolinebreak\ref{steigung}. Using this method, we are able to deduce
smoothness of the graphs and to estimate the Lipschitz constant. However,
due to the limited size of $f(U)$,
this estimate holds only on a smaller radius $\varrho<r$.
Lemma \nolinebreak \ref{iterateembedding} shows a method how to enlarge the radius, provided
the Lipschitz constant is sufficiently small. This enables us to prove the
theorem.
 \\ \\ \\
Let us come to the first statement, the comparison of two spaces
by distance bounds of finitely many points. The proof consists
of elementary geometry and is carried out here in full detail: \\[-0.2mm]
\begin{lemma} \label{steigung} Let $E \in G_{n,m}$, let
$v_1,\ldots,v_m\in E\subset \R^n$ be points on $E$ and $L\leq 1$ a constant.
If for the standard basis $\{e_1,\ldots,e_m\}$ of\,
$\R^m$ \beqn \label{abstand} |v_j-(e_j,0)|\leq
\frac{1}{3\sqrt{m}}\,L \hspace{10mm} \text{ for all }
j\in\{1,\ldots,m\}, \eeqn then $E$ is a graph over $\R^m \times\{0\}$,
that is there exists an $A=(a_1,\ldots,a_m)\in \R^{k\times m}$ with
\beqnn E=\text{\emph{span}}\{(e_1,a_1), \ldots, (e_m,a_m)\}, \eeqnn
and moreover  \beqn \label{abscha}
\|A\|=\Biggl(\sum_{j=1}^m|a_j|^2\Biggr)^{\frac{1}{2}}\leq L. \eeqn
\end{lemma}
\textbf{Proof:} \\
First we show that $E$ is a graph over
$\R^m\times \{0\}$. Suppose $E$ might not be written as a graph over
$\R^m \times \{0\}$. If $\pi$ denotes the standard projection from
$\R^n=\R^m\times \R^k$ onto $\R^m$, then \beqn \label{dimensionpie} 0 \leq \dim \pi(E) \leq m-1. \eeqn
We split the points $v_j$ into $v_j=(v_j^h,v_j^v)\in \R^m
\times \R^k$. Then, on the one hand \beqn v_1^h, \ldots, v_m^h \in
\pi(E), \eeqn and on the other hand with (\ref{abstand}) and
$L\leq 1$ for each $j$ \beqn \label{abschunterraum}
|v_j^h-e_j| &<& \frac{1}{\sqrt{m}}. \eeqn The following constructions
are carried out within the subspace $\R^m\cong\R^m\times\{0\}\subset
\R^n$. By (\ref{dimensionpie}) there exists an $e\neq 0$ in the orthogonal complement
 $[\pi(E)]^\bot \subset \R^m$. Set
$G:=\text{span}\{e\}$. Now consider the cube $Q:=[-1,1]^m\subset \R^m$ centered at the origin.
Then there is an
$s=(s_1,\ldots,s_m)\in \R^m$ with $G\cap \partial Q=\{-s,s\}$ and hence also a
$\nu \in\{1,\ldots,m\}$ with $|s_\nu|=1$. Without loss of generality
$s_\nu=1$, otherwise pass to $-s$. \\ \\ As long as $s
\neq e_\nu$ the points $0$, $e_\nu$ and $s$ constitute a rectangular triangle
with hypotenuse in $G$. The splitting
$e_\nu=e_\nu^\top + e_\nu^\bot\in G\oplus G^\bot$ yields with the
Euclidean theorem $|e_\nu|^2=|e_\nu^\top||s|$, hence \beqnn
|e_\nu^\bot-e_\nu|=\frac{1}{|s|} \geq \frac{1}{\sqrt{m}}, \eeqnn
also in the case $s=e_\nu$. But then $|w-e_\nu|\geq
\frac{1}{\sqrt{m}}$ for \emph{all} $w \in G^\bot$ and as
$\pi(E)\subset G^\bot$ \beqn |w-e_\nu|\geq
\frac{1}{\sqrt{m}} \text{ for \emph{all} } w \in \pi(E). \eeqn
But (\ref{abschunterraum}) is true also for $j=\nu$, a contradiction.
This shows that $E$ is a graph over
$\R^m\times \{0\}$. \\ \\ \\
We like to estimate the norm of $A$.
For $x\in \R^n$ and $j\in \{1,\ldots,m\}$ let $x^j\in \R^n$ be the orthogonal projection of $x$ onto
span$\{(e_j,0)\}\subset \R^n$. With $L \leq1$ and (\ref{abstand})
we have $|v_j^j- (e_j,0)| \leq |v_j- (e_j,0)| \leq
\frac{1}{3\sqrt{m}}L \leq \frac{1}{3}$, hence $|v_j^j|\geq
\frac{2}{3}$ for $1\leq j \leq m$. Let
$w_j:=\frac{1}{|v_j^j|}v_j\in E$. The second intercept theorem implies
\beqnn |w_j-w_j^j|&=&\frac{|w_j|}{|v_j|}|v_j-v_j^j| \\
&\leq& \frac{3}{2}|v_j- (e_j,0)| \\
&\leq& \frac{1}{2\sqrt{m}}L. \eeqnn With $w_j^j= (e_j,0)$ we obtain
\beqn \label{L1} |w_j- (e_j,0)|\leq \frac{1}{2\sqrt{m}}L. \eeqn Next choose
 $\nu$, such that $|a_j|\leq |a_\nu|$ for all $j$.
Without loss of generality $\nu=1$. There are
$\lambda_1,\ldots,\lambda_m \in \R$ with
$w_1=\sum_{j=1}^{m}\lambda_j(e_j,a_j)$. As $w_1^1=(e_1,0)$
we have $\lambda_1=1$. It follows \beqn \label{L2} |w_1-(e_1,0)|^2 &=&
\Biggl\lvert(0,a_1)+\sum_{j=2}^{m}\lambda_j(e_j,a_j)\Biggr\rvert^2 \nonumber \\
&=&
\Biggl\lvert\sum_{j=2}^{m}\lambda_j(e_j,0)+(0,a_1)+\sum_{j=2}^{m}\lambda_j(0,a_j)\Biggr\rvert^2
\\
&=& \sum_{j=2}^{m}\lambda_j^2+\Biggl\lvert a_1+\sum_{j=2}^{m}\lambda_ja_j\Biggr\rvert^2.
\nonumber \eeqn With (\ref{L2}), (\ref{L1}) and
$L\leq 1$ we estimate \beqn \label{L3}
\sum_{j=2}^{m}|\lambda_j|&\leq& \sqrt{m}\Biggl(\sum_{j=2}^{m}
\lambda_j^2\Biggr)^{\frac{1}{2}} \nonumber \\
&\leq& \sqrt{m}\:|w_1-(e_1,0)| \\
&\leq& \frac{1}{2}, \nonumber \eeqn and \beqn \label{L4}
\Biggl\lvert a_1+\sum_{j=2}^{m}\lambda_ja_j\Biggr\rvert\leq \frac{1}{2\sqrt{m}}L. \eeqn With
(\ref{L3}), with consideration of $|a_j|\leq |a_1|$ for all $j$,
it follows \beqn \label{L5} \Biggl\lvert\sum_{j=2}^{m}\lambda_ja_j\Biggr\rvert\leq
\Biggl(\sum_{j=2}^{m}|\lambda_j|\Biggr)|a_1| \leq \frac{1}{2}|a_1|. \eeqn From
(\ref{L5}) and (\ref{L4}) we deduce by means of absorption \beqn
\label{L6}|a_j|\leq |a_1| \leq \frac{1}{\sqrt{m}}L \hspace{10mm}
\text{for all } j \eeqn and finally $\|A\|=
\left(\sum_{j=1}^m|a_j|^2\right)^{\frac{1}{2}}\leq L$.
\hfill $\square$
\vspace{1.73cm} \\
If $f:M\rightarrow \R^n$ is an immersion and $q\in M$, then the Euclidean isometry $A_q$ is not uniquely
determined (as remarked in Section 2).
We say that a Euclidean isometry is
\emph{admissible for the point} $q\in M$, if the origin is mapped to $f(q)$
and the subspace $\R^m\times\{0\}\subset \R^m\times
\R^k$ onto $f(q)+\tau_f(q)$.
\\ \\
If a statement is true for \emph{one} admissible isometry, it often is also true for
\emph{any} admissible isometry. This will be used in the proof of the following lemma.
Although the statement of the lemma is not very surprising, its proof is quite complex as we have
to use the precise Definition \ref{definition1rlambda} in the conclusion.
Of course the numbers in the lemma are not optimal, but they suffice to
prove Theorem \ref{embedding}.
In Step 2 below we
shall apply Lemma \ref{steigung}, however the main application of this lemma
will be in the proof of Theorem \ref{embedding}.
 \\ \\
\begin{lemma} \label{iterateembedding}
Every $(r,\lambda)$-immersion with $\lambda \leq \frac{1}{8\sqrt{m}}$
is also a $(\frac{7}{4}r,8\sqrt{m}\lambda)$-immersion.
\end{lemma}
\textbf{Proof:} \\
Let $f:M^m\rightarrow \R^n$ be an $(r,\lambda)$-immersion with
$\lambda \leq \frac{1}{8 \sqrt{m}}$. \\ \\ \\
\textbf{Step 1:} Let $q \in M$, $p \in U_{r,q}$ and
$\varphi_q:=\pi\circ A_q^{-1}\circ f$, where $A_q$ is an arbitrary but fixed
admissible isometry as explained above. Then
$B_{\frac{4}{5}r}(\varphi_q(p))\subset \varphi_q(U_{r,p})$. \\ \\
Proof of Step 1: \\[1mm]
Without loss of generality we may assume
$A_q=\text{Id}_{\R^{n}}$. The set $A_q^{-1}\circ
f(U_{r,q})$ is the graph of a $C^1$-function $u:B_r\rightarrow \R^k$.
We set $w:=\varphi_q(p)\in B_r$. After a suitable rotation,
we may assume that $\{v_1,\ldots,v_m\}$ with
$v_j:=\frac{(e_j,\partial_j u(w))}{\sqrt{1+|\partial_ju(w)|^{2}}}$
for $1\leq j \leq m$
is an orthonormal basis of $\tau_f(p)$ \:(and still may assume
$A_q=\Id_{\R^n}$). Let $R \in \mathbb{SO}(n)$
be a rotation with $R(e_j,0)=v_j$ for all $j \in \{1,\ldots,m\}$.
In particular the mapping $A_p:\R^n \rightarrow \R^n$,
$A_p(x):=Rx+f(p)$, is an admissible Euclidean isometry for the point
$p \in M$. Therefore $A_p^{-1}\circ f(U_{r,p})$ is the graph of a
$C^1$-function $\tu:B_r\rightarrow \R^k$ with $\tu(0)=0$ and
$\|D\tu\|_{C^{0}(B_{r})}\leq \lambda$. We define a mapping
\beqnn g:
\overline{B}_{\frac{29}{30}r}(0)&\rightarrow&\R^m, \\
y&\mapsto& y-\pi \circ R(y,\tu(y)). \eeqnn For $y,z \in
\overline{B}_{\frac{29}{30}r}(0)$ we estimate \beqnn
|g(y)-g(z)|
&\leq& |(y-z)-\pi \circ R(y-z,0)|+|\pi \circ R(0,\tu(y)-\tu(z))| \\
&\leq& \Biggl\lvert\sum_{j=1}^{m}(y_j-z_j)(e_j-\pi v_j)\Biggr\rvert+|\tu(y)-\tu(z)| \\
&\leq& \Biggl\lvert\sum_{j=1}^{m}(y_j-z_j)
\left(1-\frac{1}{\sqrt{1+|\partial_ju(w)|^{2}}}\right)e_j\Biggr\rvert+\lambda|y-z|
\eeqnn \beqnn
&\leq&
\left(1-\frac{1}{\sqrt{1+\lambda^{2}}}\right)|y-z|+\lambda|y-z|
\\ &\leq&
(\lambda^2+\lambda)|y-z|
\\
&<& \frac{1}{6}|y-z|, \eeqnn
where we used in the last line $\lambda\leq \frac{1}{8}$.
As $g(0)=0$ we have in particular
$g(y)\in B_{\frac{1}{6}r}(0)$ for all
$y \in \overline{B}_{\frac{29}{30}r}(0)$. \\ \\
Now let $x\in \R^m$ be a point in $B_{\frac{4}{5}r}(\varphi_q(p))$.
We set $x':=x-\varphi_q(p)$. Then we have $x'
\in B_{\frac{4}{5}r}(0)$ and by the considerations above the mapping \beqnn
g+x':\overline{B}_{\frac{29}{30}r}(0)&\rightarrow&
\overline{B}_{\frac{29}{30}r}(0), \\
y &\mapsto& g(y)+x' \eeqnn is a contraction of the set
$\overline{B}_{\frac{29}{30}r}(0)$. By the Banach fixed
point theorem there is exactly one $y' \in
\overline{B}_{\frac{29}{30}r}(0)$ with $g(y')+x'=y'$, that is with $\pi
\circ R(y',\tu(y'))=x'$. Furthermore, as $y'\in B_r(0)$, there exists a
$p' \in U_{r,p}$ with $f(p')=A_p(y',\tu(y'))$. Using
$A_q=\text{Id}_{\R^{n}}$, we obtain \beqnn \varphi_q(p')
&=& \pi \circ A_q^{-1} \circ A_p(y',\tu(y')) \\
&=& \pi \circ R(y',\tu(y')) + \pi \circ A_q^{-1} \circ f(p) \\
&=& x'+\varphi_q(p) \\
&=& x. \eeqnn As $x \in B_{\frac{4}{5}r}(\varphi_q(p))$ is an arbitrary point, it follows
$B_{\frac{4}{5}r}(\varphi_q(p))\subset
\varphi_q(U_{r,p})$. \\ \\ \\[3mm]
\textbf{Step 2:} The set
$U:=U_{r,p}\cap \varphi_q^{-1}(B_{\frac{4}{5}r}(\varphi_q(p)))$ is connected
and $A_q^{-1}\circ f(U)$ is the graph of a
$C^1$-function $\hat{u}:B_{\frac{4}{5}r}(\varphi_q(p))\rightarrow \R^k$
with $\|D\hat{u}\|_{C^{0}(B_{\frac{4}{5}r}(\varphi(q)))}\leq 8
\sqrt{m}\lambda$. \\ \\
Proof of Step 2: \\[1mm]
By Step 1 we have $\pi \circ A_q^{-1}\circ
f(U)=B_{\frac{4}{5}r}(\varphi_q(p))$.
Moreover, as one can replace $\overline{B}_{\frac{29}{30}r}(0)$ in Step 1 by $\overline{B}_{r-\varepsilon}(0)$
for any sufficiently small $\varepsilon>0$, we deduce with the fixed point argument of Step 1
that $A_q^{-1}\circ f(U)$ is a graph over $B_{\frac{4}{5}r}(\varphi_q(p))$.
Now let $p'\in U_{r,p}$.
We write $A_p^{-1}\circ f(U_{r,p})$ (where $A_p$ is as in Step 1)
as graph of the $C^1$-function $\tu:B_r\rightarrow
\R^k$. Then there is a unique $x \in B_r$ with $A_p^{-1} \circ
f(p')=(x,\tu(x))$. With the rotation $R$ of Step 1 we have \beqnn
R^{-1}(\tau_f(p'))=\text{span}\{(e_1,\partial_1\tu(x)),
\ldots,(e_m,\partial_m \tu(x))\}. \eeqnn In particular \beqnn
R(e_j,\partial_j \tu(x)) &\in& \tau_f(p') \hspace{1cm} \text{ for all }
j\in \{1,\ldots,m\}. \eeqnn
Let $v_j$ and $w$ be as in Step 1. We note that $R(e_j,0)=v_j$
and estimate \beqnn |R(e_j,\partial_j \tu(x))-(e_j,0)|
&\leq& |R(e_j,\partial_j
\tu(x))-R(e_j,0)|+|R(e_j,0)-(e_j,\partial_ju(w))|\\
& & +\,|(e_j,\partial_ju(w))-(e_j,0)| \\
&=&
|\partial_j\tu(x)|+\left(\sqrt{1+|\partial_ju(w)|^2}-1\right)+|\partial_ju(w)|
\\
&\leq& 2\lambda + \left(\sqrt{1+\lambda^2}-1\right) \\
&\leq& \frac{5}{2} \lambda \\
&<& \frac{1}{3\sqrt{m}}8\sqrt{m}\lambda. \eeqnn
We apply Lemma \ref{steigung} with $E:=\tau_f(p')\in G_{n,m}$,
$v_j:=R(e_j,\partial_j\tilde{u}(x))$, $L:=8\sqrt{m}\lambda$ and
conclude that $\tau_f(p')$ may be written as a graph over $\R^m\times\{0\}$.
As this is true for any $p'\in U_{r,p}$, an argument similar to the one in the paragraph preceding (\ref{c0estforDu})
together with the considerations at the beginning of the proof of Step 2
allows us to conclude that $A_q^{-1}\circ f(U)$ is the graph of a
$C^1$-function $\hat{u}:B_{\frac{4}{5}r}(\varphi_q(p))\rightarrow \R^k$ with
$\|D\hat{u}\|_{C^{0}(B_{\frac{4}{5}r}(\varphi(q)))}\leq 8
\sqrt{m}\lambda$. In particular $\varphi_q:U\rightarrow
B_{\frac{4}{5}r}(\varphi_q(p))$ is a diffeomorphism; hence $U$ is
connected.
 \\ \\ \\
\textbf{Step 3:} The function $f$ is a
$(\frac{7}{4}r,8\sqrt{m}\lambda)$-immersion. \\ \\
Proof of Step 3: \\[1mm]
Let $\varphi_q$ and $A_q$ be as in Step 1. For every $x \in
\partial B_{\frac{19}{20}r}$ there is exactly one $p_x \in U_{r,q}$ with
$\varphi_q(p_x)=x$. For each $x \in
\partial B_{\frac{19}{20}r}$ set
$U_x:=U_{r,p_{x}}\cap \varphi_q^{-1}(B_{\frac{4}{5}r}(x))$. Moreover set
 $V_q:= U_{r,q}\cup \bigcup_{x \in \partial
B_{\frac{19}{20}r}}U_x$. By Step 1 we have
\beqnn \varphi_q(V_q) &=& B_r(0)\cup \bigcup_{x \in \partial
B_{\frac{19}{20}r}} (\varphi_q(U_{r,p_{x}})\cap B_{\frac{4}{5}r}(x))
\\
&=& B_r(0) \cup \bigcup_{x \in \partial B_{\frac{19}{20}r}}
B_{\frac{4}{5}r}(x) \\
&=& B_{\frac{7}{4}r}(0). \eeqnn
Each set $U_x$ is connected and we have
$p_x \in U_{r,q}\cap U_x$. Therefore also
$V_q$ is connected, and we have $q\in V_q$. \\ \\
Now let $R>r$ be the greatest radius, such that
$A_q^{-1}\circ f(U_{R,q})$ is the graph of a $C^1$-function
$u:B_R\rightarrow \R^k$. Suppose $R<\frac{7}{4}r$. As
$R>r$, we have $\varrho:=R-\frac{4}{5}r>0$. Define sets $U_x$
as above, but here for $x \in \partial B_\varrho$. Set
$W_q:=U_{r,q}\cup \bigcup_{x\in \partial B_\varrho}U_x$. Analogous to the
considerations above,
$W_q$ is a connected set containing $q$, and it holds
$\varphi_q(W_q)=B_R(0)$. We deduce $W_q\subset U_{R,q}$. As we assumed here,
that $A_q^{-1}\circ f(U_{R,q})$ is a graph over $B_R(0)$, and as $\varphi_q(W_q)=B_R(0)$,
we conclude $W_q=U_{R,q}$. As
$R$ is maximal, we deduce $\|Du\|_{C^{0}(B_{R})}=\infty$. But this contradicts Step 2,
saying that $\|Du\|_{C^{0}(B_{\frac{4}{5}r}(x))}\leq 8
\sqrt{m}\lambda$ for all $x \in \partial B_\varrho$.
Hence it holds $R\geq \frac{7}{4}$. \\ \\
Using the preceding considerations, we conclude
$V_q=U_{\frac{7}{4}r,q}$ (in particular $V_q$ does not depend
on the choice of $A_q$)
and $A_q^{-1}\circ f(U_{\frac{7}{4}r,q})$ is the
graph of a $C^1$-function $u:B_{\frac{7}{4}r}\rightarrow \R^k$ with
$\|Du\|_{C^{0}(B_{\frac{7}{4}r})}\leq 8 \sqrt{m}\lambda$. \\ \\[2mm]
As this is true for any point $q \in M$, the function $f$
is a $(\frac{7}{4}r,8\sqrt{m}\lambda)$-immersion.
\hfill $\square$
\vspace{0.62cm} \\
Finally we need the following lemma (which was shown in \cite{langer} for
$(r,\lambda)$-immersions): \vspace{0.2cm}
\begin{lemma} \label{inclusionc0} Let $f:M\rightarrow \R^n$ be a
$C^0$-$(r,\lambda)$-immersion and $p,q\in M$. \begin{itemize}
\item[a)] If\, $0<\varrho \leq r$ and $p \in U_{\varrho,q}$, then
$|f(q)-f(p)|<\varrho+ r\lambda$. \nopagebreak
\item[b)] If\, $\lambda\leq \frac{1}{10}$ and $p \in U_{\frac{2}{5}r,q}$,
then $U_{\frac{2}{5}r,q}\subset U_{r,p}$. \end{itemize}
\end{lemma}
\textbf{Proof:}
\begin{itemize}
\item[a)] Pass to the graph representation, use the bound on the
$C^0$-norm and the triangular inequality.
\item[b)] Let $x \in U_{\frac{2}{5}r,q}$ and $\varphi_p=\pi \circ
A_p^{-1}\circ f$. With part a) we estimate \beqnn |\varphi_p(x)|
&\leq& |f(x)-f(p)| \\
&\leq& |f(x)-f(q)|+|f(q)-f(p)| \\
&<& 2\left(\frac{2}{5}r+\frac{r}{10}\right) \\
&=& r. \eeqnn Hence $U_{\frac{2}{5}r,q}\subset
\varphi_p^{-1}(B_r)$. But $U_{\frac{2}{5}r,q}$ is a
connected set containing $p$, hence included in the $p$-component
of $\varphi_p^{-1}(B_r)$, that is in $U_{r,p}$. Hence
$U_{\frac{2}{5}r,q}\subset U_{r,p}$. \hfill $\square$
\end{itemize}
\end{section} \vspace{5mm}
\begin{section}{Proof of the embedding theorem} \vspace{1mm}
With Lemmas \ref{steigung}, \ref{iterateembedding} and \ref{inclusionc0}
we have all necessary tools for showing our theorem: \\ \\[3mm]
\noindent \textbf{Proof of Theorem \ref{embedding}:} \\
Let $m \in \N$. Define
$\Lambda=\Lambda(m):=10^{-5}m^{-2}$.
\\ \\  Now let $\lambda\leq \Lambda$,\, $r>0$ and
$f:M^m\rightarrow \R^n$ be a given
$C^{0}$-$(r,\lambda)$-immersion. We set $\varrho:=\frac{r}{5}$.
Moreover let $q \in M$ be an arbitrary point. As $2\varrho<r$, the set
$f(U_{2\varrho,q})$ may be written over
$f(q)+\tau_f(q)$ as the graph of a function
$u:B_{2\varrho}\rightarrow \R^k$ with
$\|u\|_{C^{0}(B_{2\varrho})}\leq r\lambda$. \\ \\
As the argumentation of this proof is invariant under rotations and translations,
we may assume without loss of generality that $A_q=\text{Id}_{\R^{n}}$
 (where $A_q:\R^n\rightarrow \R^n$ is an admissible isometry
for the point $q \in M$).
In particular $f(q)=0$ and
$\tau_f(q)=\R^m \times \{0\}\subset \R^m \times \R^k=\R^n$.
\\ \\ \\
Now let $x \in B_\varrho$ be an arbitrary point. Then there is exactly one
$p \in U_{\varrho,q}$ with \beqn \label{D1}
f(p)=A_q(x,u(x))=(x,u(x)). \eeqn As $\lambda \leq \frac{1}{10}$,
$2\varrho=\frac{2}{5}r$ and as $p \in U_{\varrho,q}\subset
U_{\frac{2}{5}r,q}$\,, Lemma \ref{inclusionc0} b) implies \beqnn
U_{2\varrho,q}\subset U_{r,p}. \eeqnn
Therefore the set
$f(U_{2\varrho,q})$ may be written also over $f(p)+\tau_f(p)$ as graph of a
function with small $C^0$-norm --- more precisely there exists a function
$\tu:B_r\rightarrow \R^k$ with $\|\tu\|_{C^{0}(B_{r})}\leq
r\lambda$ and $f(U_{2\varrho,q})\subset \{A_p(y,\tu(y)):y \in B_r\}$. \\ \\
Let $\{e_1,\ldots,e_m\}$ be the standard basis of $\R^m$.
For $1 \leq j \leq m$ define \beqn \label{defxj} x_j:=x+\varrho
e_j. \eeqn As $x \in B_\varrho$ we have $x_j \in B_{2\varrho}$ for
each $j$. Hence for each $j$ there is exactly one $p_j \in
U_{2\varrho,q}$ with
\beqn \label{D2} f(p_j)=A_q(x_j,u(x_j))=(x_j,u(x_j)). \eeqn \\
As $p_j \in U_{2\varrho,q}$ and $U_{2\varrho,q}\subset
 U_{r,p}$\,,
there are also unique $y_j\in B_r$ with \beqn \label{D3}
f(p_j)=A_p(y_j,\tu(y_j)). \eeqn \\ Now we estimate as follows: \beqnn
|A_p(y_j,0)-f(p)-\varrho(e_j,0)| &\leq&
|A_p(y_j,0)-f(p_j)|+|f(p_j)-f(p)-\varrho(e_j,0)| \\
&=& |A_p(y_j,0)-A_p(y_j,\tu
(y_j))|+|(x_j,u(x_j))-(x,u(x))-\varrho(e_j,0)| \\
&=& |\tu (y_j)| + |u(x_j)-u(x)| \\
&\leq& 3r\lambda \\
&=& 3\cdot 10^{-5}m^{-2}r \frac{\lambda}{\Lambda} \\
&\leq& \frac{\varrho}{3 \sqrt{m}}\cdot
8^{-3}m^{-\frac{3}{2}}\frac{\lambda}{\Lambda}. \eeqnn We
divide the inequality by $\varrho$ and obtain \beqn
\label{voraussetzungincl}
\left\lvert\frac{1}{\varrho}[A_p(y_j,0)-f(p)]-(e_j,0)\right\rvert\leq
\frac{1}{3\sqrt{m}}\cdot 8^{-3}m^{-\frac{3}{2}}\frac{\lambda}{\Lambda}. \eeqn
\\ The isometry $A_p$ maps the subspace $\R^m\times\{0\}\subset \R^m\times \R^k$
onto $f(p)+\tau_f(p)$, in particular \beqnn
\frac{1}{\varrho}[A_p(y_j,0)-f(p)] \;\in \; \tau_f(p) \hspace{8mm}
\text{ for all } j\in\{1,\ldots,m\}. \eeqnn Furthermore with $\lambda \leq
\Lambda$ we have $8^{-3}m^{-\frac{3}{2}}\frac{\lambda}{
\Lambda} \leq 1$. Hence
(\ref{voraussetzungincl}) allows us to apply Lemma \ref{steigung} with
$E:=\tau_f(p) \in G_{n,m}$,
$v_j:=\frac{1}{\varrho}[A_p(y_j,0)-f(p)]$ and
$L:=8^{-3}m^{-\frac{3}{2}}\frac{\lambda}{\Lambda}$. We conclude that
$\tau_f(p)$ may be written as a graph over $\R^m \times \{0\}$. As
$f(p)=(x,u(x))$, the implicit function theorem implies that
$u$ is differentiable in a neighborhood of $x$ and
$\tau_f(p)=\text{span}\{(e_1,\partial_1u(x)),
\ldots,(e_m,\partial_mu(x))\}$. With (\ref{abscha}) it follows \beqnn
\|Du(x)\|\leq 8^{-3} m^{-\frac{3}{2}} \frac{\lambda}{\Lambda}.
\eeqnn
\\ As $x \in B_\varrho$ was assumed to be an arbitrary point,
$u$ is differentiable on all of $B_\varrho$ and \beqn \label{c0estforDu}
\|Du\|_{C^{0}(B_{\varrho})}\leq
8^{-3}m^{-\frac{3}{2}}\frac{\lambda}{\Lambda}. \eeqn Hence, as
$\varrho=\frac{r}{5}$, the function $f$ is an
$(\frac{r}{5},8^{-3}m^{-\frac{3}{2}}\frac{\lambda}{
\Lambda})$-immersion. \\ \\
Now we can iterate the embedding
of Lemma \ref{iterateembedding} three times. Hence $f$ is also a
$(\left(\frac{7}{4}\right)^3\frac{r}{5},\frac{\lambda}{
\Lambda})$-immersion, and as $\left(\frac{7}{4}\right)^3>5$ also an
$(r,\frac{\lambda}{\Lambda})$-immersion.
This is the desired conclusion. \hfill $\square$
\end{section} \\[-2.5mm]
\begin{appendix}
\section[Graph representations over other $m$-spaces]{Graph representations over other \boldmath$m$\unboldmath-spaces}
In this appendix we would like to consider immersions with uniform graph representations
not over the affine tangent space, but over other appropriately chosen $m$-spaces.
We will show
that our theorem does not hold for such kind of immersions. \\ \\
For a given $q \in M$ and a given $m$-space $E\in G_{n,m}$ let $A_{q,E}: \mathbb{R}^n\rightarrow
\mathbb{R}^n$ be a Euclidean isometry, which maps the origin to
$f(q)$, and the subspace $\mathbb{R}^m\times\{0\}\subset
\mathbb{R}^m \times \mathbb{R}^k$ onto $f(q)+E$.
Let $U_{r,q}^E\subset M$ be the $q$-component
of the set $(\pi \circ A_{q,E}^{-1} \circ f)^{-1}(B_r)$.
Again the isometry $A_{q,E}$ is not uniquely determined but the set
$U_{r,q}^E$ does not depend on the choice of $A_{q,E}$. \\ \\
The following definition is a natural generalization of
Definition \ref{definition1rlambda}: \vspace{1mm}
\begin{defi} \label{defgeneralizedimmersion}
An immersion $f$ is called a generalized
$(r,\lambda)$-immersion, if for each point $q\in M$ there is an $E=E(q)\in G_{n,m}$, such that the set
 $A_{q,E}^{-1}\circ f(U_{r,q}^E)$ is the graph of a differentiable function
 $u:B_r\rightarrow \R^k$ with $\|Du\|_{C^{0}(B_{r})} \leq
\lambda$.
\end{defi} \vspace{1mm}
Obviously every $(r,\lambda)$-immersion is a generalized $(r,\lambda)$-immersion, as we can choose
$E(q)=\tau_f(q)$ for any $q\in M$. As a generalization of Definition \ref{c0rlambda} we have the following
definition: \vspace{1mm}
\begin{defi} \label{defi}
An immersion $f$ is called a generalized
$C^0$-$(r,\lambda)$-immersion, if for each point $q\in M$ there is an $E=E(q)\in G_{n,m}$, such that the set
 $A_{q,E}^{-1}\circ f(U_{r,q}^E)$ is the graph of a continuous function
 $u:B_r\rightarrow \R^k$ with $\|u\|_{C^{0}(B_{r})} \leq
r\lambda$.
\end{defi} \vspace{1mm}
We wonder whether there is a $\Lambda>0$, such that each generalized $C^0$-$(r,\lambda)$-immersion
with $\lambda\leq \Lambda$ is also a generalized $(r,\frac{\lambda}{\Lambda})$-immersion. The
following figure shows, that this is not the case:

\hspace{1.8cm}
\setlength{\unitlength}{0.5cm}
\begin{picture}(12,15)

\qbezier(9.2,5.1)(9.8,5.1)(9.8,5)
\qbezier(9.8,5)(9.8,4.9)(9.9,4.9)
\qbezier(9.9,4.9)(10,4.9)(10,5)

\qbezier(10,5)(10,5.1)(10.1,5.1)
\qbezier(10.1,5.1)(10.2,5.1)(10.2,5)
\qbezier(10.2,5)(10.2,4.9)(10.3,4.9)
\qbezier(10.3,4.9)(10.4,4.9)(10.4,5)

\qbezier(10.4,5)(10.4,5.1)(10.5,5.1)
\qbezier(10.5,5.1)(10.6,5.1)(10.6,5)
\qbezier(10.6,5)(10.6,4.9)(10.7,4.9)
\qbezier(10.7,4.9)(10.8,4.9)(10.8,5)

\qbezier(10.8,5)(10.8,5.1)(10.9,5.1)
\qbezier(10.9,5.1)(11,5.1)(11,5)
\qbezier(11,5)(11,4.9)(11.1,4.9)
\qbezier(11.1,4.9)(11.2,4.9)(11.2,5)

\qbezier(11.2,5)(11.2,5.1)(11.3,5.1)
\qbezier(11.3,5.1)(11.4,5.1)(11.4,5)
\qbezier(11.4,5)(11.4,4.9)(11.5,4.9)
\qbezier(11.5,4.9)(11.6,4.9)(11.6,5)

\qbezier(11.6,5)(11.6,5.1)(11.7,5.1)
\qbezier(11.7,5.1)(11.8,5.1)(11.8,5)
\qbezier(11.8,5)(11.8,4.9)(11.9,4.9)
\qbezier(11.9,4.9)(12,4.9)(12,5)

\qbezier(12,5)(12,5.1)(12.6,5.1)

\put(9.2,5.1){\line(-1,0){1}}
\put(12.6,5.1){\line(1,0){1}}

\qbezier(8.2,5.1)(4.2,5.1)(4.2,9.1)
\qbezier(4.2,9.1)(4.2,13.1)(8.2,13.1)

\qbezier(13.6,5.1)(17.6,5.1)(17.6,9.1)
\qbezier(17.6,9.1)(17.6,13.1)(13.6,13.1)
\put(8.2,13.1){\line(1,0){5.4}}

\put(9.9,2){\line(1,0){2}}
\put(4.5,6.745){\line(1,-1){1.4142}}
\put(4.18,8.1){\line(0,1){2}}

\put(1.4,6.4){\vector(1,1){2.6}}
\put(2.8,5.4){\vector(4,1){2.2}}

\put(-3.9,5.6){\small here we choose graph}
\put(-4.2,5){\small representations over the}
\put(-3.8,4.4){\small affine tangent space}

\put(11.05,2.3){\vector(1,4){0.5}}
\put(10.65,2.3){\vector(-1,4){0.5}}

\put(9.83,1.9){\tiny (}
\put(11.77,1.9){\tiny )}

\put(12.5,3.8){\small here in any point we}
\put(11.9,3.2){\small choose graph representations}
\put(12.5,2.6){\small over the horizontal line}

\end{picture}
\vspace{-0.4cm} \\ \label{graphnotsmooth}
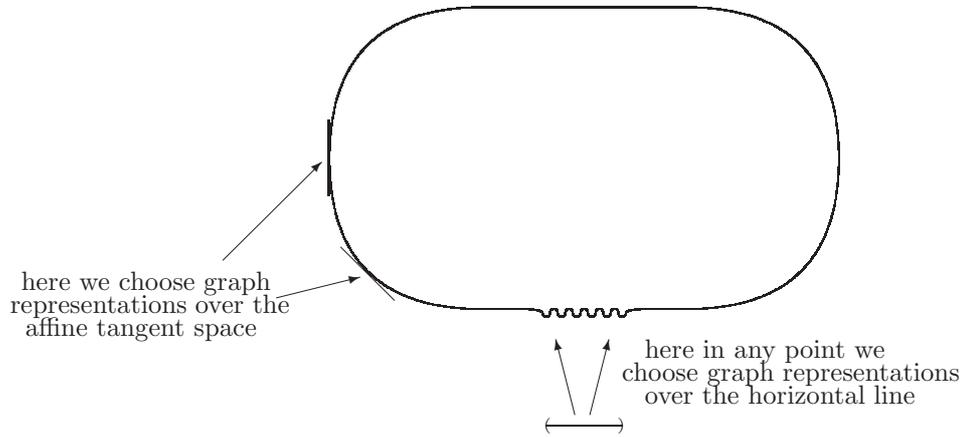
\noindent \textbf{Figure 5.1} \emph{This example shows, that a generalized $C^0$-$(r,\lambda)$-immersion
with very small $\lambda$ does not need to be a generalized $(r,\frac{\lambda}{\Lambda})$-immersion.}
\vspace{6mm} \\
Moreover, in the figure above, the part of the immersion over the horizontal line cannot be represented over any
other line (with the same radius). This shows that we require graph representations over the affine tangent space
for our theorem.
\end{appendix}

\end{document}